\documentclass[12pt]{amsart}
\usepackage{amsmath}
\usepackage{amsfonts,amssymb,amsthm, 
boxedminipage}
\usepackage{amscd}

\newtheorem{thm}{Theorem}[section]
\newtheorem{lem}[thm]{Lemma}

\newtheorem{df}[thm]{Definition}

\newtheorem{remark}[thm]{Remark}

\input epsf

\def \Z {\mathbb{Z}}

\def \CN {\mathbb N}

\begin{document}

\title[Colored Jones Polynomial]
{On the Head and the Tail of the Colored Jones Polynomial}
\author{Oliver T. Dasbach}
\address{Department of Mathematics, Louisiana State University,
Baton Rouge, LA 70803}
\email{kasten@math.lsu.edu}
\thanks {The first author was supported in part by NSF grants DMS-0306774 and DMS-0456275 (FRG)}
\author[Xiao-Song Lin]{Xiao-Song Lin}
\address{Department of Mathematics, University of California,
Riverside, CA 92521}
\thanks{The second author was supported in part by NSF grants DMS-0404511 and DMS-0456217 (FRG)}
\email{xl@math.ucr.edu}
\begin{abstract}{The colored Jones polynomial is a function 
$J_K:\CN\longrightarrow\Z[t,t^{-1}]$ associated with a knot $K$ in 3-space.
We will show that for an alternating knot $K$ the absolute values of 
the first and the last three leading coefficients of $J_K(n)$ are 
independent of $n$ 
when $n$ is sufficiently large. Computation of sample knots indicates that 
this should be true for any fixed leading coefficient 
of the colored Jones polynomial for alternating knots.
As a corollary we get a Volume-ish Theorem for the colored Jones Polynomial.
}

\end{abstract} 

\maketitle

\section {Introduction}

The celebrated Volume Conjecture of Kashaev \cite{Kashaev:VolumeConjecture} 
and Murakami and Murakami \cite{Murakami:VolumeConjecture} claims that the colored Jones polynomial 
determines the volume of a hyperbolic knot complement.
This conjecture is wide open. Recall that the colored Jones polynomial of a knot $K$ is a sequence of one variable Laurent polynomials $J_K(n)$ indexed by a positive integer $n$, the color. For $n=2$, the normalized version of 
$J_K(n)$, $J'_K(n)$, it is the classical Jones polynomial.

By using results of Lackenby, Agol and Thurston \cite{Lackenby:Volume} we showed in
\cite{DL:VolumeIsh} that there is a Volume-ish Theorem for the Jones polynomial of 
alternating knots: The volume of a hyperbolic alternating knot is linearity bounded by
the sum of the absolute values of the second leading and the second lowest coefficient
of the Jones polynomial $J'_K(2)$.

Here we will extend our results to the colored Jones polynomial of alternating knots.
We will show that in this case the colored Jones polynomial has a well defined head and
tail: The leading three and the lowest three coefficients are independent of the color, for $n \geq 3$.
It turns out that these coefficients have a particular nice form and can
easily be computed. It is interesting to note, though, that the third coefficient of the
colored Jones polynomial $J_K'(3)$ is in general not determined by the Jones polynomial $J_K'(2)$.

Moreover, the coefficients that gave rise to the Volume-ish Theorem for the Jones polynomial
persists in all the colored Jones polynomials: The second leading coefficient and the second lowest coefficient are independent of the color.
As a corollary we will get a 
Volume-ish Theorem for the colored Jones polynomial, independent of the color $n$.
Thus there is indeed a deeper connection between the colored Jones polynomial 
and the volume of hyperbolic knots, at least for alternating knots.

\noindent
{\bf Acknowledgment: } 
We are more than indebted to Dror Bar-Natan for making his wonderful Mathematica package KnotTheory publicly available on his web page \cite{Bar-Natan:KnotTheory}.
The colored Jones polynomial routine allowed us to experimentally observe some of the phenomena that we were able to show here.

\section{Background}

\begin{figure}[ht] 
\centerline{\epsfxsize=3in\epsfbox{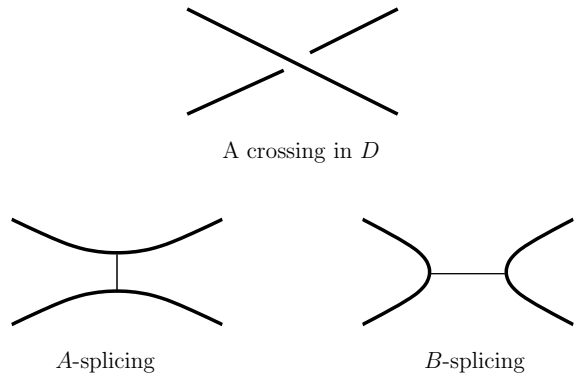}}
\caption{Splicings of a crossing, $A$-graph and $B$-graph.\label {AB-splicing}}
\end{figure}

Let $D$ be a knot diagram. For each crossing of $D$, we may apply two kinds of
splicings (see Figure \ref{AB-splicing}) as in the Kauffman skein relation. 
One splicing is called $A$-splicing and the other $B$-splicing. A state $s$
of $D$ is a choice of $A$-splicing or $B$-splicing for each crossing of 
$D$. After applying a state $s$ to splice each crossing of $D$, we change $D$
to a diagram of a collection of disjoint simple closed curves (circles)
in the plane. Let
$d(s)$ be the number of the resulting circles. Also let $\alpha(s)$ 
($\beta(s)$, resp.) the number of $A$-splicings ($B$-splicings, resp.) 
in the state $s$.
With these notations, 
the Kauffman bracket of $D$ is
$$\langle D\rangle = \sum_s\,A^{\alpha(s)-\beta(s)}(-A^2-A^{-2})^{d(s)-1},$$
where the sum is over all states $s$ of $D$.

\medskip
\medskip

\begin{df} The {\rm $A$-graph} $A(D)$ of a knot diagram $D$ has its set 
of vertices equal to the set of circles in the all $A$-splicing state, and
there is an edge
joining two vertices of $A(D)$ for each crossing in $D$ 
between the corresponding circles. Similarly, we define the {\rm $B$-graph} 
$B(D)$ of $D$ using the all $B$-splicing state. 
\end{df}
 
To the best of our knowledge, the notion of $A$-graph and $B$-graph of a knot diagram $D$ was 
first introduced by A. Stoimenow (\cite{Stoimenow:SecondCoefficient}, compare with \cite{Thistlethwaite:AdequateKauffman}).

\begin{df} A knot diagram $D$ is called $A$-adequate ($B$-adequate, resp.)
if there are no
loops, i.e. edges that start and terminate at the same vertex, in the $A$-graph $A(D)$ (the $B$-graph $B(D)$, resp.). 
A knot diagram $D$ is called adequate if it is both $A$-adequate and 
$B$-adequate. 
\end{df}

It is easy to see that a reduced alternating knot diagram is adequate. 
In fact, in this case the $A$- and the $B$-graph are the two checkerboard graphs
of the knot diagram.

\begin{remark}
A diagram $D$ of a knot is $A$-adequate if and only if the corresponding diagram $D^*$ of its mirror image is $B$-adequate.
Furthermore, the Kauffman bracket of the mirror image $D^*$ equals the Kauffman bracket of $D$ after substituting $A^*$ for the variable $A$. 
Therefore, all theorems on $A$-adequate knots have a corresponding
formulation  for $B$-adequate knots. In the sequel we will frequently omit this formulation to make the theorems more
readable. 
\end{remark}

The reduced $A$-graph $A(D)'$ of a knot diagram is obtained from
the $A$-graph $A(D)$ by keeping the same set of vertices but reduce 
all multiple edges to one for each pair of vertices. So a knot diagram $D$ 
is $A$-adequate if and only if $A(D)'$ is a simple graph, i.e. does not have loops or multiple edges. 
Similarly, we have the reduced
$B$-graph $B(D)'$ and $D$ is $B$-adequate if and only if $B(D)'$ is a simple graph.

Let $D$ is a knot diagram of $c=c(D)$ crossings. 
The number of vertices of $A(D)'$, which equals the number of vertices in the unreduced graph $A(D)$, is denoted by $v=v(D)$ and the number of edges of $A(D)'$ by $e=e(D)$. 

The following theorem is a generalization of our corresponding result in \cite{DL:VolumeIsh} that we showed for alternating knots:

\begin{thm} {\rm (Stoimenow)}\label{stoim} 
Suppose $D$ is $A$-adequate. Then we have
$$\langle D\rangle = (-1)^{v-1}A^{c+2v-2}+(-1)^{v-2}(e-v+1)A^{c+2v-6}+\text{\rm
lower order terms}.$$
\end{thm}

\begin{proof} Recall that $e=e(D)$ is the number of edges of the reduced 
graph $A(D)'$. Let $k_i$ be the multiplicity of the $i$-th edge of $A(D)'$
in $A(D)$, for $i=1,2,\dots,e$.

For a state $s$ of $D$ where there are at least two $B$-splicings appearing
at edges of $A(D)$ which reduce to different edges in $A(D)'$, its contribution
to $\langle D\rangle$ is 
$$A^{c-2\beta(s)}(-A^2-A^{-2})^{d(s)-1}=(-1)^{d(s)-1}
A^{c-2\beta(s)+2d(s)-2}+\text{
lower degree terms},$$ 
where $d(s)\leq v-2+\beta(s)-2$. Now
$$c-2\beta(s)+2d(s)-2\leq c+2v-10.$$

So we only need to consider states whose $B$-splicings all appear 
on the same $k_i$ multiple edges of $A(D)$. The contribution of such 
states in $\langle D\rangle$ is
$$\begin{aligned}
&A^c(-A^2-A^{-2})^{v-1}+\sum_{i=1}^e\sum_{j=1}^{k_i}\,
\binom{k_i}{j} A^{c-2j}\,\,
(-A^2-A^{-2})^{v-2+j-1}\\
&=(-1)^{v-1}\left(A^{c+2v-2}+(v-1)A^{c+2v-6}+\cdots\right)\\
&\quad+\sum_{i=1}^e\sum_{j=1}^{k_i}\,\binom{k_i}{j}\,(-1)^{v+j-3}A^{c+2v-6}+
\cdots\\
&=(-1)^{v-1}A^{c+2v-2}+(-1)^{v-2}(e-v+1)A^{c+2v-6}+\text{lower degree terms}.
\end{aligned}$$  
This proves the theorem.
\end{proof}

Note that when the graph $A(D)'$ is connected, we have
$$e-v+1=\beta_1(A(D)') \quad\text{
(the first Betti number of $A(D)'$).}$$ This will be important to us.

\begin{lem}\label{cab-b-1} 
Let $D$ be an $A$-adequate (resp. $B$-adequate) knot diagram and $D^n$ be 
its blackboard $n$ cabling. Then $D^n$ is $A$-adequate (resp. $B$-adequate). 
Furthermore, we have
$$\beta_1(A(D^n)')=\beta_1(A(D)').$$
\end{lem}

\begin{proof} The first claim was shown by Lickorish (see e.g. \cite{Lickorish:Book}). The following
picture depicts the $A$-splicing of the blackboard cabling of a knot diagram
$D$:

\medskip
\begin{figure}[ht]
\centerline{\epsfxsize=3in\epsfbox{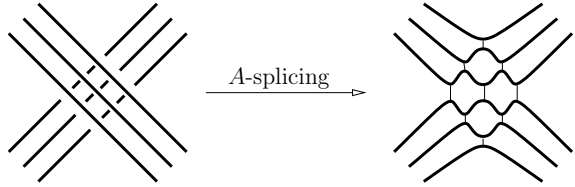}}
\caption{$A$-splicing of the cabling of a crossing.}
\end{figure}
\medskip

If $D$ is $A$-adequate, then $A(D)'$ is a simple graph. For a vertex $V$ of 
$A(D)'$, let the edges coming out of $V$ be $E_1,E_2,\dots,E_k$. Then in 
$A(D^n)'$, we have $n$ vertices associated with $V$: $V_1,V_2,\dots,V_n$, 
and $n-1$ edges $V_1V_2,V_2V_3,\dots, V_{n-1}V_n$. The edges 
$E_1,E_2,\dots,E_k$ of $A(D)'$ are still in $A(D^n)'$ as edges coming out of 
$V_1$. Thus, $A(D^n)'$ is equal $A(D)'$ with a chain 
$V_1V_2\cup V_2V_3\cup\cdots\cup V_{n-1}V_n$ attached to each vertex $V=V_1$.
So $A(D^n)'$ is still a simple graph and $D^n$ is $A$-adequate.

The second claim follows easily from the structure of $A(D^n)'$ described
above.
\end{proof}

\section{The colored Jones polynomial}

We will use the Chebyshev basis of the Kauffman bracket skein module of
$S^1\times[0,1]$ to express the colored Jones polynomial. See e.g. 
\cite{Lickorish:Book}.

Let $S_n(x)$, $n\geq0$, be polynomials of $x$ specified by the following 
recurrence relation and initial values:
$$S_{n+1}=xS_n-S_{n-1},\quad S_0(x)=1,\quad S_1(x) = x.$$
Inductively, we have
\begin{equation}\label{cheby}
S_n(x)=x^n+(1-n)x^{n-2}+\text{lower degree terms}.
\end{equation}

Suppose that $D$ is a knot digram of a knot $K$. 
Let $J_K(n)$ be the colored Jones polynomial of $K$ corresponding to the 
$n$ dimensional irreducible $U_q(\mathfrak{sl}_2)$ module. 
Here, $J_K(2)$ is the classical Jones polynomial multiplied by $(A^2+A^{-2})$. We have
$$J_K(n+1)=\left[(-1)^nA^{n^2+2n}\right]^{-w(D)}\,(-1)^{n-1}\,[2]\,
\langle S_n(D)\rangle,$$
where $[2]=A^2+A^{-2}$, $w(D)$ is the writhe of $D$, and $S_n(D)$ is a
linear combination of blackboard cablings of $D$ obtained using the
Chebyshev polynomial $S_n(x)$. By Equation (\ref{cheby}), we have
\begin{equation}\label{cheby2}
S_n(D)=D^n+(1-n)D^{n-2}+\text{lower degree cablings of $D$.}
\end{equation}

For the unknot diagram $O$, we have
$$J_O(n)=\frac{A^{2n}-A^{-2n}}{A^2-A^{-2}}=A^{2n-2}+A^{2n-6}+\cdots
+A^{-2n+6}+A^{-2n+2}:=[n].$$
The normalized colored Jones polynomial is
$$J'_K(n)=J_K(n)/[n]$$ and
$J'_K(2)$ is the classical Jones polynomial.

We are interested in the leading coefficients of $J'_D(n)$.
We write

$$J'_D(n)=aA^k+bA^{k-4}+cA^{k-8}+\cdots,$$
thus
$$\begin{aligned}
J_D(n)&=J'_D(n)\,[n]\\
&=(aA^k+bA^{k-4}+cA^{k-8}+\cdots)(A^{2n-2}+A^{2n-6}+\cdots)\\
&=aA^{k+2n-2}+(a+b)A^{k+2n-6}+(a+b+c)A^{k+2n-10}+\cdots
\end{aligned}
$$

Suppose now that $D$ is an $A$-adequate knot diagram.
Since all blackboard cablings of $D$ are adequate, using 
Theorem \ref{stoim}, Lemma \ref{cab-b-1}, and Equation (\ref{cheby2}), we have
$$
\begin{aligned}
\langle S_n(D)\rangle&=\langle D^n\rangle + (1-n)\langle D^{n-2}\rangle+\cdots\\
&=(-1)^{nv-1}A^{n^2c+2nv-2}+(-1)^{nv-2}\beta_1(A(D)')A^{n^2c+2nv-6}+\cdots\\
&\quad+(-1)^{(n-2)v-1}(1-n)A^{(n-2)^2c+2(n-2)v-2}\\
&\quad+(-1)^{(n-2)v-2}(1-n)
\beta_1(A(D)')A^{(n-2)^2c+2(n-2)v-6}+\cdots
\end{aligned}$$
Compare the second highest degree of $\langle D^n\rangle$ and the highest degree
of $\langle D^{n-2}\rangle$:
$$(n^2c+2nv-6)-((n-2)^2c+2(n-2)v-2)=(4n-4)c+4v-4>0.$$
So,
$$\langle S_n(D)\rangle
=(-1)^{nv-1}A^{n^2c+2nv-2}+(-1)^{nv-2}\beta_1(A(D)')A^{n^2c+2nv-6}+\text{lower 
degree terms}.$$
Finally, we have
$$
\begin{aligned}
\,[2]\,\langle S_n(D)\rangle&=(A^2+A^{-2})\langle S_n(D)\rangle\\
&=(-1)^{nv-1}A^{n^2c+2nv}+((-1)^{nv-1}+(-1)^{nv-2}\beta_1(A(D)'))A^{n^2c+2nv-4}\\
&\quad+
\text{lower 
degree terms}.
\end{aligned}
$$

\begin{thm} Let $D$ be an $A$-adequate knot diagram. Write
$$J'_D(n)= a_n A^{k_n}+b_nA^{k_n-4}+\text{\rm lower degree terms}.$$
Then, we have
$$|a_n|=1\quad\text{and}\quad|b_n|=\beta_1(A(D)').$$
\end{thm}

\begin{proof} This follows directly from the previous calculations.
\end{proof}

\noindent{\bf Remark:} Note that the difference between highest degrees of 
$\langle D^n\rangle$ and $\langle D^{n-2}\rangle$, respectively, 
is $O(n)$. So when $n$ is getting large, any fixed leading portion of 
coefficients of $J'_K(n)$ will only depend on the coefficients of $\langle 
D^n\rangle$. For the first and second leading coefficients of 
$\langle D^n\rangle$, we have shown that they depend only on $n$ by a sign
when $K$ has an $A$-adequate diagram.
We consider the third leading coefficient of $\langle D^n\rangle$ in the 
next section.

\section{The third coefficient of $J'_K(n)$} \label{Section Third Coefficient}

The difficulties in computing the third leading coefficient are that it does not 
depend on the topological type of the $A$-graph and its reduced graph anymore. 
Informations on the embedding of the $A$-graph are needed.

We first will extend Theorem \ref{stoim}. Suppose that a knot diagram
$D$ is $A$-adequate. The vertices of $A(D)$ are represented by disjoint 
circles in the plane. Let $E$ and $E'$ be two different 
edges in $A(D)'$. We say 
that they are {\it disjoint} if $E\cap E'=\emptyset$. If $E$ and $E'$ share 
a common vertex $V$, then there are two possibilities. In one case, the
circle representing $V$ can be divided into two semi-circles such that 
vertices of multiple edges $E$ and $E'$ lie on different semi-circles 
respectively. We say that $E$ and $E'$ are {\it separated} at $V$ in this 
case. Otherwise, we say that $E$ and $E'$ are {\it mixed} at $V$. See Figure 3,
in which the pairs $\{E_1,E_4\}$ are $\{E_1,E_5\}$ are mixed, and the
pairs $\{E_4,E_5\}$ and $\{E_4,E_6\}$ are separated. 

\medskip
\begin{figure}[ht]
\centerline{\epsfxsize=3in\epsfbox{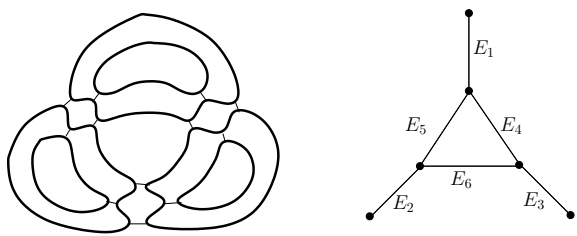}}
\caption{$A(D)$ and $A(D)'$ where $D$ is the 2-cable of the
left-trefoil diagram.}
\end{figure}
\medskip

Suppose $E$ and $E'$ are two edges of $A(D)'$ mixed at a vertex $V$. Then in 
$A(D)$, the multiple edge $E$ ($E'$, resp.) is partitioned into sets of 
$p_1,p_2,\dots,p_m$ ($q_1,q_2,\dots,q_m$, resp.) parallel multiple edges,
so that the endpoints of these edges on the circle $V$ are placed 
alternatively as $p_1,q_1,p_2,q_2,\dots,p_m,q_m$. We must have $m\geq2$. 
See Figure 4. If 
$D$ is a reduced alternating knot diagram, we do not have such pairs
of edges $E$ and $E'$ in $A(D)'$. 

\medskip
\begin{figure}[ht]
\label{mixededges}
\centerline{\epsfxsize=3in\epsfbox{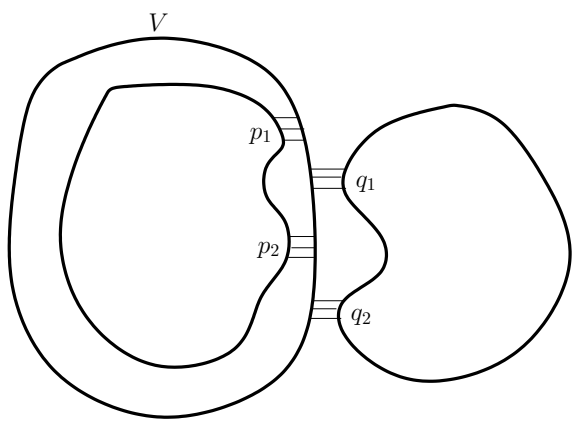}}
\caption{Two multiple edges mixed at $V$.}
\end{figure}
\medskip

\begin{thm} \label{Theorem41}
Suppose $D$ is $A$-adequate. Then we have
$$\begin{aligned}
\langle D\rangle &= (-1)^{v-1}A^{c+2v-2}+(-1)^{v-2}(e-v+1)A^{c+2v-6}\\
&\quad+(-1)^{v-3}\left(\binom{v-1}{2}-e(v-2)+\mu+\binom{e}{2}-\theta-
\tau\right)A^{c+2v-10}\\
&\quad+\text{\rm
lower order terms,}
\end{aligned}$$ 
where $\mu$ is the number of edges in $A(D)'$ whose multiplicity in $A(D)$
is larger than 1, $\theta$ is the number of pairs of edges in $A(D)'$ mixed at 
a vertex, and $\tau$ is the number of triangles in $A(D)'$.
\end{thm}

When $D$ is a reduced alternating knot diagram, then $\theta=0$. So this theorem 
generalized the formula of $\langle D\rangle$ in \cite{DL:VolumeIsh} for reduced 
alternating knot diagrams.

\begin{proof} 

Let $s$ be a state of $D$ with at least three $B$-splicings such that
(1) they appear at edges of $A(D)$ which reduce to three different edges 
in $A(D)'$, and (2) these three edges of $A(D)'$ 
do not form a triangle in $A(D)'$. 
Then the contribution of $s$
to $\langle D\rangle$ is 
$$A^{c-2\beta(s)}(-A^2-A^{-2})^{d(s)-1}=(-1)^{d(s)-1}A^{c-2\beta(s)+2d(s)-2}+\text{
lower degree terms},$$ 
where $d(s)\leq v-3+\beta(s)-3$. Now
$$c-2\beta(s)+2d(s)-2\leq c+2v-14.$$

So we only need to consider the following four cases.

\medskip
\noindent 
{\bf Case 1:} States whose $B$-splicings all appear 
on the same $k_i$ multiple edges of $A(D)$. The total contribution of such 
states in $\langle D\rangle$ is
$$\begin{aligned}
&A^c(-A^2-A^{-2})^{v-1}+\sum_{i=1}^e\sum_{j=1}^{k_i}\,
\binom{k_i}{j} A^{c-2j}\,\,
(-A^2-A^{-2})^{v-2+j-1}\\
&=(-1)^{v-1}\left(A^{c+2v-2}+(v-1)A^{c+2v-6}+\binom{v-1}{2}A^{c+2v-10}+
\cdots\right)\\
&\quad+\sum_{i=1}^e\sum_{j=1}^{k_i}\,\binom{k_i}{j}\,(-1)^{v+j-3}A^{c+2v-6}+
\sum_{i=1}^e\sum_{j=1}^{k_i}\,\binom{k_i}{j}\,(-1)^{v+j-3}(v+j-3)A^{c+2v-10}+
\cdots\\
&=(-1)^{v-1}A^{c+2v-2}+(-1)^{v-2}(e-v+1)A^{c+2v-6}\\
&\quad +(-1)^{v-3}\left(\binom{v-1}{2}-e(v-2)+\mu\right)A^{c+2v-10}+
\text{lower degree terms}.
\end{aligned}$$  

The last equation follows from:
$$\begin{aligned}
\sum_{j=1}^{k_i} {k_i \choose j} (-1)^j j &=& \left \{ 
\begin {array}{rl}
-1 & k_i=1\\
0 & k_i>1
\end{array}
\right .
\end{aligned}
$$

\medskip
\noindent {\bf Case 2:} States whose $B$-splicings appear in a pair of 
distinct multiple edges, which are either disjoint or separated.

The total contribution of such 
states in $\langle D\rangle$ is
$$
\begin{aligned}
&\sum_{i=1}^{k_r}\sum_{j=1}^{k_s}\,
\binom{k_r}{i}\binom{k_s}{j}\,A^{c-2(i+j)}(-A^2-A^{-2})^{v-4+i+j-1}\\
&=(-1)^{v-3}A^{c+2v-10}+\text{lower order terms}.
\end{aligned}
$$

\noindent {\bf Case 3:} States whose $B$-splicings appear in a pair of 
distinct multiple edges $E,E'$ which are mixed at a vertex. This is the 
most delicate of the four cases. We claim that there is no contribution
to the three leading coefficients that involve $B$-splicings at both
edges.

Suppose the $A$-diagram is locally similar to Figure \ref{mixededges},
i.e. $$D=D(p_1, q_1, p_2, q_2, \dots ,p_k, q_k)$$ with
$p_i>0$ and $q_i>0$ for  $i=1,\dots,k$.

We will deal with the following subcases:
\begin{enumerate}
\item \label {subcase1} One of the $p_i$ or one of the $q_i$ is greater than $1$.
Assume that $p_1>1$. In this case one can see that
the three leading coefficients of $\langle D(p_1,q_1, \dots, p_k, q_k) \rangle$ equal the
three leading coefficients of  $\langle D(p_1-1,q_1, \dots, p_k, q_k) \rangle$.
Therefore, to show the claim we can assume that $p_1=\dots=p_k=q_1=\dots=q_k=1$.
\item \label {subcase2} If $k=2$, i.e. our diagram is $D(1,1,1,1)$, a straightforward computation shows that
the contribution of the type of $B$-splicings we consider here vanishes.
\item We are left with the case of $k>2$. We compare the leading terms of
$$\langle D(p_1=1, q_1=1, p_2=1, q_2=1, \dots, p_k=1, q_k=1)\rangle$$
and 
$$\begin{aligned}
& \langle D(p_1=0, q_1=1, p_2=2, q_2=1, \dots, p_k=1, q_k=1)\rangle = \\
& \langle D(p_2=2, q_2=1, \dots, p_k=1, q_k=2) \rangle.
\end{aligned}$$

Yet another direct, straightforward computation shows that the three leading
coefficients coincide. Thus we can reduce this case to (\ref{subcase1}) and
(\ref{subcase2}).
\end{enumerate}

Thus we showed that the sum of all states with $B$-splicings at a pair of mixed edges does not
contribute to the leading three coefficients.
Recall that in our notation $\theta$ is the number of pairs of edges in $A(D)'$
that are mixed. 
To summarize Case 2 and Case 3: The total contribution of states whose
$B$-splicings appear in a pair of distinct multiple edges is
$$(-1)^{v-3}\left(\binom{e}{2}-\theta\right)A^{c+2v-10}+\text{lower order terms}.
$$

\noindent {\bf Case 4:}
States with three $B$-splicings appearing on
three different edges in $A(D)'$ that form a triangle. 
Note, that for reasons of planarity the edges cannot be mixed at any of the
three vertices.

The total contribution of such states is
$$
\begin{aligned}
\sum_{\tau}&\sum_{i=1}^{k_r}\sum_{j=1}^{k_s}\sum_{k=1}^{k_t}\,\binom
{k_r}{i}\binom{k_s}{j}\binom{k_t}{k}A^{c-2(i+j+k)}
(-A^2-A^{-2})^{v-1+i+j+k-3-1}\\
&=(-1)^{v-3}(-\tau)A^{c+2v-10}+\text{lower order terms}.
\end{aligned}
$$

This proves the theorem.
\end{proof}

Similar to the discussion in \S3, we need to compare first the third highest 
degree of $\langle D^n\rangle$ and the highest degree of 
$\langle D^{n-2}\rangle$:
$$(n^2c+2nv-10)-((n-2)^2c+2(n-2)v-2)=(4n-4)c+4v-8>0$$
since we have $v\geq 2$ and $n\geq2$ for non-trivial cases. Thus, we just need 
to consider the third leading coefficient of $\langle D^n\rangle$.

Denote by $v_n$ the number of vertices of $A(D^n)'$, $e_n$ the number of 
edges in $A(D^n)'$, $\mu_n$ the number of edges in $A(D^n)$ whose 
multiplicity in $A(D^n)$ is larger than 1, $\theta_n$ the number of 
pairs of edges in $A(D^n)'$ mixed at a vertex,
and $\tau_n$ the number of 
triangles in $A(D^n)'$.

\begin{lem} Suppose the $A$-adequate
diagram $D$ contains no kinks, or $A(D)$ has no vertices of valence 1.
We have:
$v_n=n v$, $e_n=e+(n-1)v$, $\mu_n=e_n=e+(n-1)v$ when $n>1$, and $\tau_n=\tau$.
Furthermore, if $D$ is a reduced alternating knot diagram and $n>1$,
then $\theta_n=(n-2)v+2e$.
\end{lem}

\begin{proof} This is clear from the structure of $A(D^n)'$.
\end{proof}

The relevant terms in Theorem (\ref{Theorem41}) for the third coefficient are now simplified with the following Lemma:

\begin{lem} If $D$ is a reduced alternating knot diagram and $n>1$,
we have
$$\binom{v_n-1}{2}-e_n(v_n-2)+\mu_n+\binom{e_n}{2}-\theta_n-
\tau_n=\frac{(e-v)^2+(e-v)}{2}-\tau+1.$$
\end{lem}

Thus, we have finally arrived at our main theorem.

\begin{thm} \label{Main Theorem} Let $K$ be an alternating knot. Write
$$\begin{aligned}
J'_K(n)&=& \pm \left ( a_n A^{k_n}-b_nA^{k_n-4}+c_nA^{k_n-8}\right) \pm  \dots \\
& &\pm \left ( \gamma_n A^{k_n-4 r_n+8} -
\beta_n A^{k_n-4 r_n+4}+\alpha_n A^{k_n-4 r_n} \right )
\end{aligned}$$
with positive $a_n$ and $\alpha_n$.

Let $A(D)$ and $B(D)$ be the $A$- and $B$-graphs of a reduced alternating diagram
$D$ of $K$ with crossing number $c$. The reduced graphs $A(D)'$ and $B(D)'$ have $e_A$ and $e_B$ edges and $v_A$ and $v_B$ vertices. Note, that $v_A+v_B=c+2$. 
Furthermore, there are $\tau_A$ and $\tau_B$ triangles
in $A(D)'$ and $B(D)'$. 

Then
\begin{enumerate}
\item \label{span} The span $r_n$ of $J'_K(n)$, i.e. the difference of the highest and lowest exponent,  in the variable $q=A^4$ is:
$$r_n={n \choose 2} c.$$
\item $a_n=\alpha_n=1.$
\item $b_n=e_A-v_A+1$ and
$\beta_n=e_B-v_B+1$.
\item $c_n= {b_n \choose 2}-\tau_A$ and $\gamma_n={\beta_n \choose 2}-\tau_B$
for $n >2 $. The coefficients $c_2$ and $\gamma_2$ were identified in \cite{DL:VolumeIsh}:
\begin{eqnarray*}
c_2&=&{b_2+1 \choose 2} + n(2) - \tau_A\\
\gamma_2&=& {\beta_2 +1 \choose 2} + n^*(2) -\tau_B,
\end{eqnarray*}
\end{enumerate}
where $n(2)$ (or $n^*(2)$)  is the number of edges in $A'(D)$ (or $B'(D)$) that correspond to edges of multiplicity at least $2$ in
the unreduced graph $A(D)$ (or $B(D)$).

In particular, $|a_n|,|b_n|, |\alpha_n|$ and $|\beta_n|$ are independent of $n$ when $n>1$ and $|c_n|$ and $|\gamma_n|$
are independent of $n$ when $n>2$.
\end{thm}

\begin{proof}
Item (\ref{span}) was essentially shown by Kurpita and Murasugi in \cite{KurpitaMurasugi:HierarchalJones} and  Le in \cite{Le:AJConjecture}.
The highest degree of $\langle S(D^n) \rangle$ in $A$ is $n^2 c + 2 n v_a-2$ and
the lowest degree is $-n^2 c - 2 n v_b+2$.
Thus for the span of $\langle S(D^n) \rangle$ in the variable $A$ we
get
$$
2 n^2 c + 2 n (v_A+v_B) + 4
= 2 n^2 c + 2 n (c+2) +4.
$$
This implies that the span of $J_K(n+1)$ is
$$2 n^2 c + 2 n (c+2)$$
and thus the span of $J_K'(n+1)$ in the variable $A$ is
$$2 n^2 c + 2 n (c+2) - 4n = 2 n^2 c + 2 n c.$$
Changing  to the variable $q$ this shows Claim (\ref{span}).

For the remaining claims: 
Recall that we showed  
$$\begin{aligned}
(-1)^{p_n} A^{l_n}\langle S(D^n)\rangle &=
1-(e_A-v_A+1)A^{-4} \\
&\quad+ \left(\frac{(e_A-v_A)^2+(e_A-v_A)}{2}-\tau_A+1 \right) A^{-8} \\
&\quad+\text{lower order terms}
\end{aligned}
$$
for some number $p_n$ and $l_n$.

The remaining claims now follows from a comparison of the coefficients of
$\langle S(D^n) \rangle$ and $J'_K(n+1)$ and the fact that the $A$-graph of the knot is the
$B$-graph of the mirror image of the knot and vice versa.
\end{proof}

\section {Example}

\begin{figure}[ht]
\label{Figure 12_217}
\centerline{\epsfxsize=3in\epsfbox{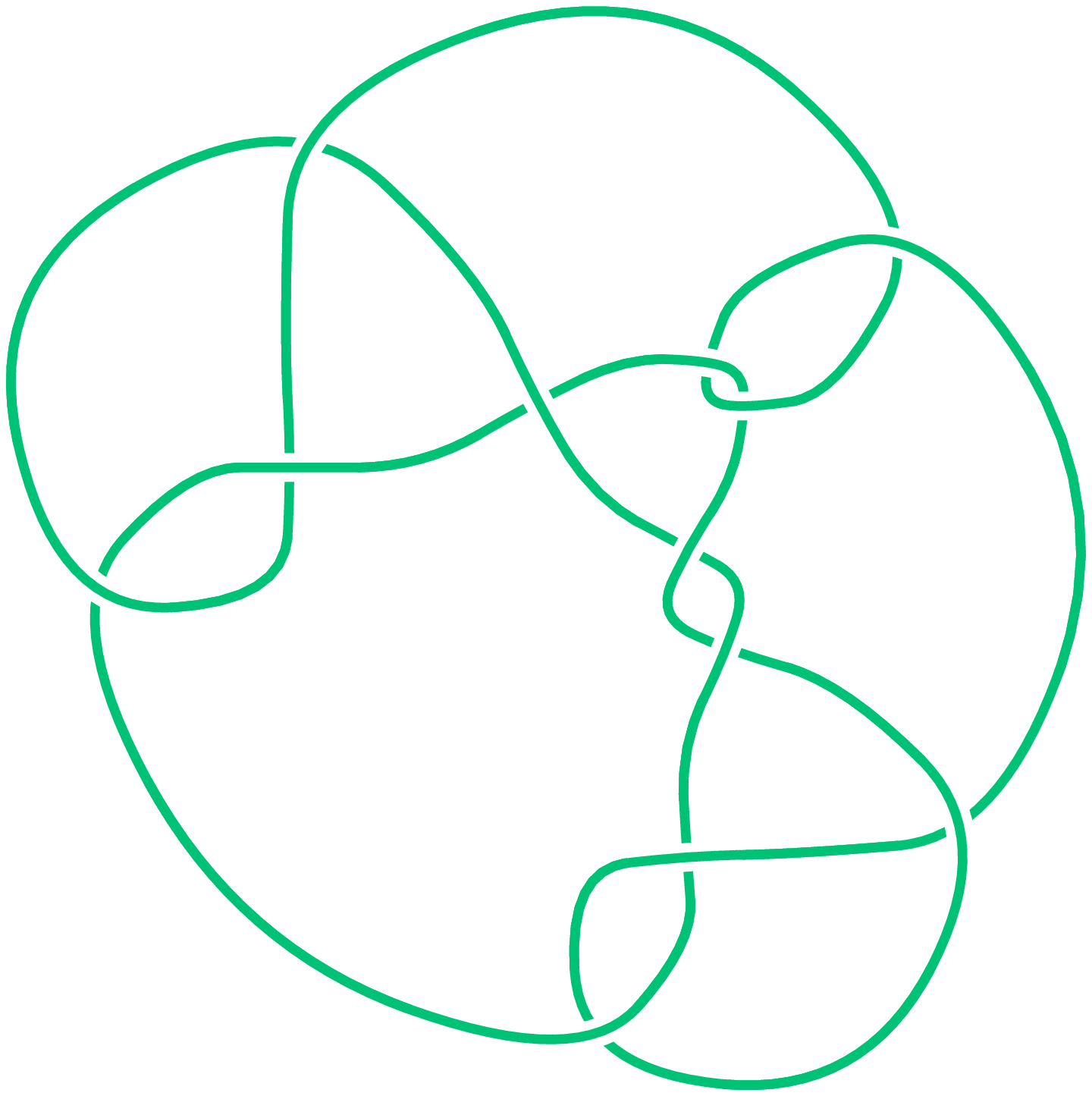}  \epsfxsize=3in\epsfbox{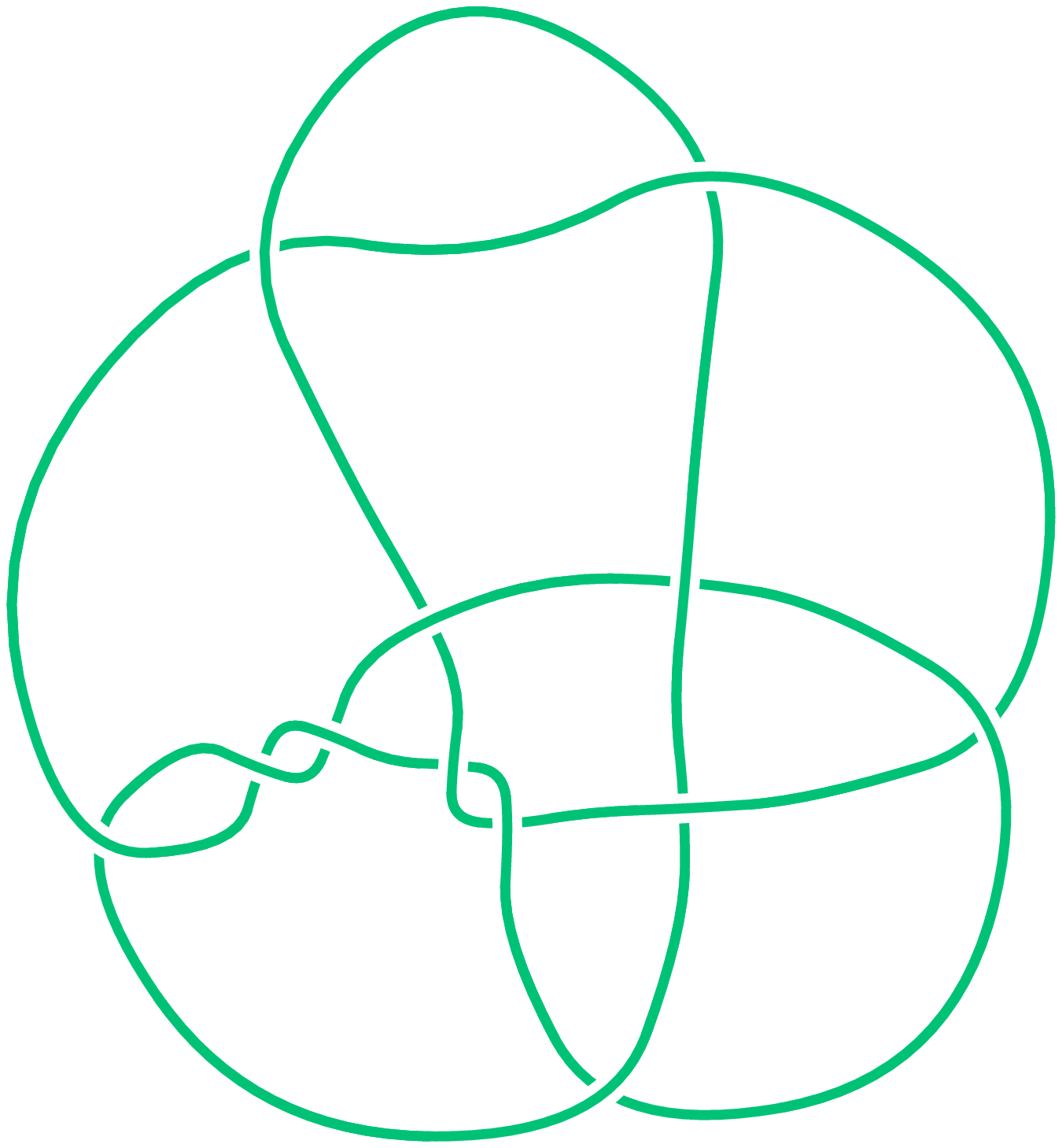} }
\caption{The alternating knots $12_{217}$ and $12_{1228}$}
\end{figure}

As an example we pick the alternating knots $12_{217}$ and $12_{1228}$ in the Knotscape census \cite{HTW:KnotTabulation}.
The two knots are given in Figure \ref{Figure 12_217}. According to Knotscape both knots share the same Jones Polynomial $J'_K(2)$.

We have for the knot $12_{217}$:

The $A$-graph is the checkerboard graph corresponding to the color not containing the outer face.
Its reduced graph has $e_A=9$ edges and $v_A=6$ vertices. Furthermore, it contains $\tau_A=4$ triangles.
Thus the three leading coefficients are, up to a common sign change: $1, -(e_A-v_A+1)=-4$ and ${4 \choose 2}- \tau_A=2.$

Similarly: $e_B=11, v_B=8$ and $\tau_B=2$. Thus the last three coefficients are: $4, -4, 1$, up to a common sign.

For the knot $12_{1228}$ the data are: $e_A=9, v_A=6, \tau_A=3, e_B= 11, v_B=8, \tau_B=2$ and
we get for the leading coefficients $1, -4, 3$ and for the last three coefficients: $4, -4, 1$.
Again, these two lists are up to sign changes. 
 
Since the Jones polynomials of the two knots coincide but already the third coefficient of the colored Jones polynomials for color $n>2$ 
does not, we see:

\begin{lem}
For $n>2$ the leading third coefficient of the colored Jones polynomial $J'_K(n)$ is in general not
determined by the Jones polynomial $J'_K(2)$.
\end{lem}

\section {The Volume-ish Theorem for the colored Jones polynomial}

As a corollary to the Main Theorem \ref{Main Theorem} and the results in \cite{DL:VolumeIsh}   we also get the Volume-ish Theorem for the
colored Jones polynomial:

\begin{thm}[Volume-ish Theorem for the colored Jones polynomial] For an alternating, prime, non-torus knot $K$ let
$$J'_K(n)=a_n q^{k_n} + b_n q^{k_n-1}+ \dots + \beta_n q^{k_n-r_n+1} + \alpha_n q^{k_n-r_n}$$
be the colored Jones polynomial.

Then with $b:=b_2$ and $\beta:=\beta_2$ we have $b=b_n$ and $\beta=\beta_n$ for all $n$ and:

$$2 v_0(max(|b|,|\beta|)-1) \leq Vol(S^3-K) \leq 10 v_0(|b|+|\beta|-1)$$
where $v_0 \approx 1.0149416$ is the volume of an ideal regular hyperbolic tetrahedron.
\end{thm}

exi
\providecommand{\bysame}{\leavevmode\hbox to3em{\hrulefill}\thinspace}
\providecommand{\MR}{\relax\ifhmode\unskip\space\fi MR }
\providecommand{\MRhref}[2]{%
  \href{http://www.ams.org/mathscinet-getitem?mr=#1}{#2}
}
\providecommand{\href}[2]{#2}

\end{document}